\documentclass[12pt,thmsa]{article}
\usepackage{amssymb}

\usepackage{sw20lart}

\input{tcilatex}
\begin{document}

\author{Alan Horwitz \\
Penn State University\\
25 Yearsley Mill Rd.\\
Media, PA 19063\\
alh4@psu.edu}
\title{Invariant Means}
\date{7/14/00}
\maketitle

\begin{abstract}
Let $m(a,b)$ and $M(a,b,c)$ be symmetric means. We say that $M$ is type 1
invariant with respect to $m$ if $M(m(a,c),m(a,b),m(b,c))\equiv M(a,b,c)$.
If $m$ is strict and isotone, then we show that there exists a unique $M$
which is type 1 invariant with respect to $m$. In particular we discuss the
invariant logarithmic mean $L_{3},$ which is type 1 invariant with respect
to $L(a,b)=\frac{b-a}{\log b-\log a}.$ We say that $M$ is type 2 invariant
with respect to $m$ if $M(a,b,m(a,b))\equiv m(a,b).$ We also prove existence
and uniqueness results for type 2 invariance, given the mean $M(a,b,c)$. The
arithmetic, geometric, and harmonic means in two and three variables satisfy
both type 1 and type 2 invariance. There are means $m$ and $M$ such that $M$
is type 2 invariant with respect to $m,$ but not type 1 invariant with
respect to $m($for example, the Lehmer means). $L_{3}$ is type 1 invariant
with respect to $L,$ but not type 2 invariant with respect to $L.$
\end{abstract}

\section{Introduction}

Let $R_{+}^{n}=\{(a_{1},...,a_{n})\in R^{n}:a_{i}>0\;\forall i\}.$ A \textit{%
mean} $m(a_{1},...,a_{n})$ in $n$ variables is a continuous function on $%
R_{+}^{n}$ with $\min (a_{1},...,a_{n})\leq m(a_{1},...,a_{n})\leq \max
(a_{1},...,a_{n})$. $m$ is called \textit{symmetric} if $m(\pi
(a_{1},...,a_{n}))=m(a_{1},...,a_{n})$ for any permutation $\pi $. A mean $%
m(a_{1},...,a_{n})$ is called:

$\cdot $ Strict if $m(a_{1},...,a_{n})=\min (a_{1},...,a_{n})$ or $%
m(a_{1},...,a_{n})=\max (a_{1},...,a_{n})$ if and only if $a_{1}=\cdots
=a_{n}$.

$\cdot $ Homogeneous if $m(ka_{1},...,ka_{n})=km(a_{1},...,a_{n})$ for any $%
k>0$

$\cdot $ Isotone(strictly) if $m(a_{1},...,a_{n})$ is an
increasing(strictly) function of each of its variables.

We let $\Sigma _{n}$ denote the set of means in $n$ variables.

Let $A_{2}(a,b)=\dfrac{a+b}{2}$ and $A_{3}(a,b,c)=\dfrac{a+b+c}{3}$ denote
the arithmetic means in two and three variables, respectively. A simple
computation shows that

$A_{3}(A_{2}(a,c),A_{2}(a,b),A_{2}(b,c))\equiv A_{3}(a,b,c).$ The same type
of invariance holds for the geometric and harmonic means--i.e., $%
G_{3}(G_{2}(a,c),G_{2}(a,b),G_{2}(b,c))\equiv G_{3}(a,b,c)$ and $%
H_{3}(H_{2}(a,c),H_{2}(a,b),H_{2}(b,c))\equiv H_{3}(a,b,c)$. The work for
this paper started with the aim of finding a mean $M(a,b,c)$ which has the
same invariance property with respect to the logarithmic mean in two
variables given by

\[
L(a,b)=\frac{b-a}{\log b-\log a}
\]

We call this mean the invariant logarithmic mean $L_{3}(a,b,c).$ In general, 
$M(a,b,c)$ is said to be invariant with respect to $m(a,b)$ if $%
M(m(a,c),m(a,b),m(b,c))\equiv M(a,b,c).$ Of course this notion can be
extended to means in any number of variables, though we concentrate in this
paper on means in two and three variables. We also discuss extensions of
some of the other classical means in two variables to invariant means in
three or more variables. It is not obvious that such invariant means even
exist for a given $m(a,b)$. We prove(Theorem \ref{T1}) that if $m$ is strict
and isotone, then an invariant mean $M(a,b,c)$ does exist. In particular,
there is an invariant logarithmic mean in three variables. Various authors
have extended the logarithmic mean to three or more variables(see \cite{H}, 
\cite{PT}, and \cite{ST}), but none of those means are invariant. We shall
compare the invariant logarithmic mean to those means.

We also show(Proposition \ref{P2}) that $M$ inherits properties of $m$ such
as symmetry, homogeneity, and isotonicity. We say that a symmetric mean $%
M(a,b,c)$ is invariant if there exists a symmetric $m(a,b)$ such that $M$ is
type 1 invariant with respect to $m$. Not every symmetric mean $M$ is
invariant. For example, we show that $M(a,b,c)=\left( \dfrac{ab+ac+bc}{3}%
\right) ^{1/2}$ is not invariant.

As noted above, the arithmetic, geometric, and harmonic means in three
variables are invariant , respectively, with respect to the arithmetic,
geometric, and harmonic means in two variables . However, consider the
well-known Lehmer means $lh_{p}(a,b)=\dfrac{a^{p}+b^{p}}{a^{p-1}+b^{p-1}}${}
and $LH_{p}(a,b,c)=\dfrac{a^{p}+b^{p}+c^{p}}{a^{p-1}+b^{p-1}+c^{p-1}}$. It
is easy to show that $LH_{p}$ is \textit{not} invariant with respect to $%
lh_{p}$ if $p\neq 1$. However, $LH_{p}$ is invariant with respect to $lh_{p}$
in the sense that $LH_{p}(a,b,lh_{p}(a,b))\equiv lh_{p}(a,b)$ . We call this
type 2 invariance. The extension of this type of invariance to means in $n$
variables is $m_{n}(a_{1},a_{2},...,a_{n-1},m_{n-1}(a_{1},...,a_{n-1}))%
\equiv m_{n-1}(a_{1},...,a_{n})$. Let us call the invariance discussed
earlier(given by (\ref{inv})) type 1 invariance. Since the Lehmer means
satisfy type 2 invariance, type 2 invariance does not imply type 1
invariance. However, type 1 invariance does not imply type 2 invariance
either. For example, $L_{3}(a,b,c)$ is type 1 invariant with respect to $%
L(a,b)$, but not type 2 invariant with respect to $L(a,b).$

Given $M(a,b,c)$, we also prove existence and uniqueness results for type 2
invariance(see Theorem \ref{e2} and Theorem \ref{uniq}). We can also prove
that given $m(a,b),$ there exists a symmetric mean $M(a,b,c)$ such that $M$
is type 2 invariant with respect to $m$. Such an $M$ is not unique, however.
We \textit{cannot} yet show that if $m$ is analytic, then $M$ can be chosen
to be analytic as well.

While the ideas in this paper can also be discussed for non-symmetric means,
we usually restrict our results to symmetric means.

\section{Preliminary Material}

In this section we give some elementary results on means and symmetric
functions which will be useful in later sections.

\begin{lemma}
\label{sym}

(1) Let $f(x,y)$ be a differentiable, \textit{symmetric} function on an open
region $E\subset R^{2}$, and assume that $D=\{(x,y)\in E:x=y\}$ is nonempty.
Let $i,j$ be nonnegative integers with $i+j=n$. Then$\dfrac{\partial ^{n}f}{%
\partial x^{i}\partial y^{j}}(a,a)=\dfrac{\partial ^{n}f}{\partial
x^{j}\partial y^{i}}(a,a)$ for any $(a,a)\in D$.

(2) Let $f(x,y,z)$ be a differentiable, \textit{symmetric} function on an
open region $E\subset R^{3}$, and assume that $D=\{(x,y,z)\in E:x=y=z\}$ is
nonempty. Let $i,j,k$ be nonnegative integers with $i+j+k=n$, and let $%
(r,s,t)$ be any permutation of $(i,j,k)$. Then$\dfrac{\partial ^{n}f}{%
\partial x^{i}\partial y^{j}\partial z^{k}}(a,a,a)=\dfrac{\partial ^{n}f}{%
\partial x^{r}\partial y^{s}\partial z^{t}}(a,a,a)$ for any $(a,a,a)\in D$.
\end{lemma}

\begin{remark}
If $Q=(a,a),$ then Lemma \ref{sym}, part (2) implies, for $n=2$, that $%
M_{xy}(Q)=M_{xz}(Q)=M_{yz}(Q)$ and $M_{xx}(Q)=M_{yy}(Q)=M_{zz}(Q)$. If $%
Q=(a,a,a)$, then for $n=3$ we have $M_{xxy}(Q)=M_{xxz}(Q)=M_{xyy}(Q)=\cdots $
and $M_{xxx}(Q)=M_{yyy}(Q)=M_{zzz}(Q)$.
\end{remark}

\begin{proposition}
\label{0.0}Let $m(a,b)$ be a differentiable, symmetric mean, and let $%
Q=(a,a),$ $a>0$. Then

(i) $m_{x}(Q)=m_{y}(Q)=\dfrac{1}{2}$

(ii) $m_{xy}(Q)=-m_{xx}(Q)$

(iii) $3m_{xxy}(Q)=-m_{xxx}(Q)$
\end{proposition}

\proof%
The proof follows by taking successively differentiating both sides of $%
m(x,x)=x.$ The details are similar to the proof of Proposition \ref{0.1}
below, and we omit them.

\begin{theorem}
\label{m}Let $m(a,b)$ be a differentiable, symmetric, and homogeneous mean,
and let $f(x)=m(a,x),$ $a>0$. Then

(i) $f^{\prime }(a)=\dfrac{1}{2}$

(ii) $f^{\prime \prime \prime }(a)=-\dfrac{3}{2a}f^{\prime \prime }(a)$
\end{theorem}

\proof%
(i) follows immediately from Proposition \ref{0.0}, \#1. To prove (ii), note
first that, since $m$ is homogeneous, $xm_{x}+ym_{y}=m$. Taking $\dfrac{%
\partial }{\partial x}$ of both sides yields 
\begin{equation}
xm_{xx}+ym_{yx}=0  \label{2.5}
\end{equation}
$\Rightarrow xm_{xx}(x,a)+am_{yx}(x,a)=0\Rightarrow m_{yx}(x,a)=-\dfrac{x}{a}%
m_{xx}(x,a)=-\dfrac{x}{a}f^{\prime \prime }(x)$

$\Rightarrow m_{yxx}(x,a)=-\dfrac{1}{a}(xf^{\prime \prime \prime
}(x)+f^{\prime \prime }(x))\Rightarrow m_{yxx}(a,a)=-\dfrac{1}{a}(af^{\prime
\prime \prime }(a)+f^{\prime \prime }(a).$ By Proposition \ref{0.0} and
Lemma \ref{sym}, $m_{yxx}(a,a)=-\dfrac{1}{3}f^{\prime \prime \prime }(a)$.
Setting $-\dfrac{1}{a}(af^{\prime \prime \prime }(a)+f^{\prime \prime }(a)=-%
\dfrac{1}{3}f^{\prime \prime \prime }(a)$ proves (ii).

\begin{proposition}
\label{0.1}Let $M(a,b,c)$ be a differentiable, symmetric mean, and let $%
Q=(a,a,a),$ $a>0$. Then

(i) $M_{x}(Q)=M_{y}(Q)=M_{z}(Q)=\dfrac{1}{3}$

(ii) $M_{xy}(Q)=-\dfrac{1}{2}M_{xx}(Q)$

(iii) $M_{xyz}(Q)=-\dfrac{1}{2}(M_{xxx}(Q)+6M_{xxy}(Q))$

(iv) $M_{xxxx}(Q)+8M_{xxxy}(Q)+6M_{xxyy}(Q)+12M_{xxyz}(Q)=0$
\end{proposition}

To prove that $M_{x}(Q)=\dfrac{1}{3}$, take $\dfrac{d}{dx}$ of both sides of
the identity $M(x,x,x)=x.$ This gives 
\begin{equation}
(M_{x}+M_{y}+M_{z})(x,x,x)=1  \label{3}
\end{equation}
(i) now follows from Lemma \ref{sym}. Taking $\dfrac{d}{dx}$ of both sides
of (\ref{3}) and using Clairut's Theorem gives 
\begin{equation}
(M_{xx}+M_{yy}+M_{zz}+2M_{xy}+2M_{xz}+2M_{yz})(x,x,x)=0  \label{3.5}
\end{equation}

Lemma \ref{sym} yields $3M_{xx}(Q)+6M_{xy}(Q)=0$, which proves (ii)$.$
Taking $\dfrac{d}{dx}$ of both sides of (\ref{3.5}) and using Clairut's
Theorem again gives\newline
$\left(
M_{xxx}+M_{yyy}+M_{zzz}+3(M_{xxy}+M_{xyy}+M_{xxz}+M_{xzz}+M_{yyz}+M_{yzz})+6M_{xyz}\right) (Q)=0
$

Lemma \ref{sym} then yields (iii). Finally, (iv) follows in a similar
fashion.

\begin{proposition}
\label{0.2}Let $M(a,b,c)$ be a differentiable, symmetric, and homogeneous
mean, and let $Q=(a,a,a),$ $a>0$. Then

(i) $M_{xx}(Q)=-a(M_{xxx}(Q)+2M_{xxy}(Q))$

(ii) $M_{xxx}(Q)=-\dfrac{a}{2}(M_{xxxx}(Q)+2M_{xxxy}(Q))$

(iii) $M_{xxy}(Q)=-\dfrac{a}{2}(M_{xxyy}(Q)+M_{xxxy}(Q)+M_{xxyz}(Q))$

(iv) $M_{xyz}(Q)=-\dfrac{3a}{2}M_{xxyz}(Q)$
\end{proposition}

\proof%
Since $M$ is homogeneous, $xM_{x}+yM_{y}+zM_{z}=M$. Taking $\dfrac{\partial 
}{\partial x}$ of both sides yields 
\begin{equation}
xM_{xx}+yM_{yx}+zM_{zx}=0  \label{4}
\end{equation}
Taking $\dfrac{\partial }{\partial x}$ of both sides of (\ref{4}) gives 
\begin{equation}
xM_{xxx}+M_{xx}+yM_{yxx}+zM_{zxx}=0  \label{4.5}
\end{equation}
Lemma \ref{sym} then implies $M_{xx}(Q)=-a(M_{xxx}+2M_{xxy})(Q),$ which
proves (i). Taking $\dfrac{\partial }{\partial x}$ or $\dfrac{\partial }{%
\partial y}$ of both sides of (\ref{4.5}) and using Lemma \ref{sym} gives
(ii) and (iii). Finally, taking $\dfrac{\partial }{\partial y}$ of both
sides of (\ref{4}), and then taking $\dfrac{\partial }{\partial z}$ of both
sides yields $xM_{xxyz}+yM_{xyyz}+zM_{xyzz}+2M_{xyz}=0$. Lemma \ref{sym}
then implies $3aM_{xxyz}+2M_{xyz}=0$, which is (iv).

\section{Type 1 Invariance\label{s2}}

\begin{definition}
A symmetric mean $m_{n}(a_{1},...,a_{n})$ is said to be \textit{invariant}
with respect to $m_{n-1}(a_{1},...,a_{n-1})$ if 
\begin{equation}
m_{n}(m_{n-1}(a_{2},...,a_{n}),m_{n-1}(a_{1},a_{3},...,a_{n}),...,m_{n-1}(a_{1},...,a_{n-1}))=m_{n}(a_{1},...,a_{n})
\label{inv}
\end{equation}
for all $(a_{1},...,a_{n})\in R_{+}^{n}$.\footnote{%
Later we call this type 1 invariance.}
\end{definition}

For example, if $A_{n}(a_{1},...,a_{n})=\dfrac{a_{1}+\cdots +a_{n}}{n}$ and $%
G_{n}(a_{1},...,a_{n})=(a_{1}\cdots a_{n})^{1/n}$ denote the arithmetic and
geometric means, respectively, in $n$ variables, then for each $n\geq 3$,%
\newline
$%
A_{n}(A_{n-1}(a_{2},...,a_{n}),A_{n-1}(a_{1},a_{3},...,a_{n}),...,A_{n-1}(a_{1},...,a_{n-1}))\equiv A_{n}(a_{1},...,a_{n}) 
$

and

$%
G_{n}(G_{n-1}(a_{2},...,a_{n}),G_{n-1}(a_{1},a_{3},...,a_{n}),...,G_{n-1}(a_{1},...,a_{n-1}))\equiv G_{n}(a_{1},...,a_{n}). 
$

For the rest of this section we discuss means $M(a,b,c)$ which are invariant
with respect to a given mean $m(a,b)$. In that case, we have 
\begin{equation}
M(m(a,c),m(a,b),m(b,c))=M(a,b,c)  \label{inv3}
\end{equation}
Also, given $M(a,b,c)$, if a mean $m(a,b)$ exists with $M$ invariant for $m$%
, we say that $M$ is an \textit{invariant mean}. We shall have more to say
about this ``reverse'' process in \S \ref{reverse}.

In this section, we are given a mean $m(a,b)$, and we assume throughout that 
$m$ is a strict, isotone, and symmetric mean. The strictness and isotonicity
are necessary in general in order for our proofs to work.

First, given positive numbers $a,b,c$, define the following recursive
sequences:

\begin{equation}
a_{n+1}=m(a_{n},c_{n}),c_{n+1}=m(a_{n},b_{n}),b_{n+1}=m(c_{n},b_{n}),a_{0}=a,b_{0}=b,c_{0}=c
\label{1}
\end{equation}

We now prove that the sequences $\{a_{n}\},\{b_{n}\},\{c_{n}\}$ each
converge to the same limit, which lies in the smallest interval containing $%
a,b,$ and $c$.

\begin{proposition}
\label{P1} Each of the sequences $\{a_{n}\},\{b_{n}\},\{c_{n}\}$ converges,
and

$\lim_{n\rightarrow \infty }a_{n}=\lim_{n\rightarrow \infty
}b_{n}=\lim_{n\rightarrow \infty }c_{n}=L$, with $\min
(a_{0},b_{0},c_{0})\leq L\leq \max (a_{0},b_{0},c_{0})$.
\end{proposition}

\proof%
Assume, without loss of generality, that $a_{0}\leq c_{0}\leq b_{0}$. We
claim

\begin{equation}
a_{k}\leq c_{k}\leq b_{k}\;\forall k  \label{2}
\end{equation}

We prove (\ref{2}) by induction. So assume that $a_{n}\leq c_{n}\leq b_{n}$
for some $n\geq 0$. Then $a_{n+1}=m(a_{n},c_{n})\leq m(a_{n},b_{n})=c_{n+1}$
and $b_{n+1}=m(c_{n},b_{n})\geq m(a_{n},b_{n})=c_{n+1}$. That proves (\ref{2}%
).

Now for each $n$, $a_{n+1}\geq a_{n}$ since $a_{n}\leq c_{n}$, and $%
b_{n+1}\leq b_{n}$ since $c_{n}\leq b_{n}$. Using (\ref{2}), this implies
that $a_{n}\leq b_{n}\leq b_{0}$ and $b_{n}\geq a_{n}\geq a_{0}$. Since $%
\{a_{n}\}$ is increasing and bounded above, and $\{b_{n}\}$ is decreasing
and bounded below, $\{a_{n}\}$ converges to $L_{1}\geq a_{0}$ and $\{b_{n}\}$
converges to $L_{2}\leq b_{0}$. Since $c_{n+1}=m(a_{n},b_{n}),$ $\{c_{n}\}$
converges to $m(L_{1},L_{2})\equiv L_{3}$. Now $a_{n+1}=m(a_{n},c_{n})%
\Rightarrow L_{1}=m(L_{1},L_{3})\Rightarrow L_{1}=L_{3}$ since $m$ is
strict. Similarly, $L_{2}=L_{3}$. Letting $L$ equal the common value of $%
L_{1}$,$L_{2}$, and $L_{3},$ each of the sequences $\{a_{n}\},\{b_{n}\},%
\{c_{n}\}$ converges to $L$, with $a_{0}\leq L\leq b_{0}$.%
\endproof%

Proposition \ref{P1} shows that the recursion (\ref{1}) defines a mean in
three variables $M(a,b,c)=L$, where $a_{0}=a,b_{0}=b,c_{0}=c$, except for
the \textit{continuity} of $M$. We shall prove that in Theorem \ref{cont}
below. Henceforth we use $M(a,b,c)$ to denote the common limit of the
sequences defined by (\ref{1}). Of course $M$ depends on the given mean $m$.
We prove next that an invariant mean for $m$ is defined precisely in this
way.

\begin{theorem}
\label{T1}Let $L$ denote the common limit of the sequences defined by (\ref
{1}).

(A) Define the mean, $M$, in three variables by $M(a,b,c)=L$, where $%
a_{0}=a,b_{0}=b,c_{0}=c$. Then $M$ is an invariant mean for $m$.

(B) If $M$ is an invariant mean for $m$, then $M(a,b,c)=L$, where $%
a_{0}=a,b_{0}=b,c_{0}=c$.
\end{theorem}

\proof%
(A) $M(m(a,c),m(a,b),m(c,b))=$ common limit of the sequences defined by (\ref
{1}), with $a_{0}=m(a,c),b_{0}=m(a,b),c_{0}=m(c,b)$. But this is the same
limit as that of the sequences defined by (\ref{1}), with $%
a_{0}=a,b_{0}=b,c_{0}=c$. Hence $M(m(a,c),m(a,b),m(c,b))=M(a,b,c)$.

(B) Since $M$ is an invariant mean for $m$, $%
M(a_{0},b_{0},c_{0})=M(a_{1},b_{1},c_{1})=\cdots =M(a_{n},b_{n},c_{n})$.
Taking the limit as $n\rightarrow \infty $ gives $%
M(a_{0},b_{0},c_{0})=M(L,L,L)=L$. 
\endproof%

Since the limit of a sequence is unique, Theorem \ref{T1} yields the
following.

\begin{corollary}
\label{C1}Let $M_{1}$ and $M_{2}$ be invariant means for $m$. Then $%
M_{1}=M_{2}$.
\end{corollary}

In light of the corollary, we can now speak of \textit{the }invariant mean
for $m$.

\begin{remark}
Our approach above is similar in many ways to the well known idea of
compounding three given means $M_{1},M_{2},M_{3}$ in three variables to
obtain another mean $[M_{1},M_{2},M_{3}]$ in three variables(see \cite{B}).
Indeed, the invariant mean $M$ can be obtained by compounding the means $%
M_{1}(a,b,c)=m(a,c),$ $M_{2}(a,b,c)=m(a,b),\;$and $M_{3}(a,b,c)=m(b,c)$.
However, the standard theorems on compound means do not appear to imply
Proposition \ref{P1} or Theorem \ref{T1}. For means in two variables, the
existence of the compound mean $[M_{1},M_{2}]$ is proved, for example, in 
\cite{B} with the assumption that $M_{1}$ and $M_{2}$ are \textit{comparable}%
. A similar assumption is left out of the theorem in \cite{B} for the
existence of $[M_{1},M_{2},M_{3}].$ However, the proof given also seems to
require the assumption of comparability. Note, however, that the means $%
M_{1}(a,b,c)=m(a,c),$ $M_{2}(a,b,c)=m(a,b),\;$and $M_{3}(a,b,c)=m(b,c)$ are 
\textbf{not comparable }in general.\textbf{\ }.
\end{remark}

Next we show that $M$ inherits many properties of $m$.

\begin{proposition}
\label{P2}Suppose that $M(a,b,c)$ is invariant for $m(a,b)$.

(A) Then $M$ is strict, isotone, and symmetric.

(B) If $m$ is homogeneous, then $M$ is homogeneous.

(C) Suppose that $m_{1}(a,b)\leq m(a,b)\leq m_{2}(a,b)$ for all $(a,b)\in
R_{+}^{2}$, and let $M_{1}$ and $M_{2}$ be the invariant means for $m_{1}$
and $m_{2}$ respectively. Then $M_{1}(a,b,c)\leq M(a,b,c)\leq M_{2}(a,b,c)$
for all $(a,b,c)\in R_{+}^{3}$.
\end{proposition}

\proof%
The isotonicity of $M$ follows immediately from the isotonicity of $m$, as
does the homogeneity of $M$ if $m$ is homogeneous. To prove the symmetry of $%
M$ we use (\ref{1}). If we have the permutation $a_{0}\leftrightarrow b_{0}$%
, then by (\ref{1}) and the symmetry of $m$, $a_{n}\leftrightarrow b_{n}$
for all $n$, and hence the common limit $L$ of the three sequences remains
the same. Thus $M(a,b,c)=M(b,a,c)$. If we have the permutation $%
a_{0}\leftrightarrow c_{0}$, then again by (\ref{1}) and the symmetry of $m$%
, $c_{1}\leftrightarrow b_{1}$, $c_{2}\leftrightarrow a_{2}$, $%
c_{3}\leftrightarrow b_{3},...$ Hence the sequence $\{a_{n}\}$ gets sent to
the sequence $c_{0},a_{1},c_{2},a_{3},c_{4},a_{5},...$, which still
converges to $L$. Thus $M(c,b,a)=M(a,b,c)$. Similarly, one can show that $%
M(a,c,b)=M(a,b,c)$. Hence $M(\pi (a,b,c))=M(a,b,c)$ for any permutation $\pi 
$.

Now we prove that $M$ is strict. So suppose that $a_{0}<c_{0}\leq b_{0}$.
Then $a_{1}=m(a_{0},c_{0})>a_{0}$ since $m(a,b)$ is strict. This implies
that $L>a_{0}$ since $\{a_{n}\}$ is increasing by the proof of Proposition 
\ref{P1}. Now $c_{1}=m(a_{0},b_{0})\Rightarrow a_{0}<c_{1}<b_{0}$ since $%
m(a,b)$ is strict. Then $b_{1}=m(c_{0},b_{0})\leq b_{0}$ and $%
b_{2}=m(c_{1},b_{1})\leq m(c_{1},b_{0})<b_{0}$, which implies that $L<b_{0}$
since $\{b_{n}\}$ is decreasing, again by the proof of Proposition \ref{P1}.
Thus we have proven (A) and (B).

To prove (C), let $\{a_{n}^{(k)}\},\{b_{n}^{(k)}\},\{c_{n}^{(k)}\}$ denote
the sequences defined by (\ref{1}), with $m=m_{k},\;k=1,2$, and starting
with the same initial values $a_{0}=a,b_{0}=b,c_{0}=c$. It is not hard to
show, using induction, that $a_{n}^{(1)}\leq a_{n}\leq
a_{n}^{(2)},\;b_{n}^{(1)}\leq b_{n}\leq b_{n}^{(2)},\;$and $c_{n}^{(1)}\leq
c_{n}\leq c_{n}^{(2)}$. It then follows immediately that $M_{1}(a,b,c)\leq
M(a,b,c)\leq M_{2}(a,b,c).$

\endproof%

Given $m(a,b)$, we now define the map $\phi :$ $\sum {}_{3}\rightarrow $ $%
\sum {}_{3}$

\[
\phi (N)(a,b,c)=N(m(a,c),m(a,b),m(c,b)),\;N\in \sum {}_{3}
\]

It follows immediately that a mean $M\in \sum {}_{3}$ is a fixed point of $%
\phi $ if and only if $M$ is invariant with respect to $m$. In addition, we
now prove that the iterates $\phi ^{[n]}(N)$ converge to the invariant mean $%
M$.

\begin{theorem}
\label{T2}For any $N\in \sum {}_{3}$ ,\ $\lim_{n\rightarrow \infty }\phi
^{[n]}(N)$ exists and equals the invariant mean for $m$.
\end{theorem}

\proof%
$\phi ^{[n]}(N)(a,b,c)=N(a_{n},b_{n},c_{n})$, where $a_{0}=a,b_{0}=b,c_{0}=c$%
, and the sequences $\{a_{n}\},\{b_{n}\},\{c_{n}\}$ are defined by (\ref{1}%
). By Proposition \ref{P1}, $\lim_{n\rightarrow \infty }\phi
^{[n]}(N)(a,b,c)=N(L,L,L)=L$ for each $(a,b,c)\in R_{+}^{3}$, where $L$ is
the common limit of the three sequences. The Theorem now follows from
Theorem \ref{T1}, Part A. 
\endproof%

The following result shall prove useful when comparing invariant means to
other known means.

\begin{theorem}
\label{T4}Let $M$ be the invariant mean for $m$, and let $N(a,b,c)$ be any
mean in threee variables. If $N(m(a,c),m(a,b),m(c,b))\leq (\geq )N(a,b,c)$
for all $(a,b,c)\in R_{+}^{3}$, then $M(a,b,c)\leq (\geq )N(a,b,c)$ for all $%
(a,b,c)\in R_{+}^{3}$. Furthermore, the sequence of means $\phi
^{[n]}(N)(a,b,c)$ is decreasing(increasing).
\end{theorem}

\proof%
We prove the $\leq $ case. For any positive integer $n$, if $\phi
(N)(a,b,c)= $\newline
$N(m(a,c),m(a,b),m(c,b))\leq N(a,b,c)$ for all $(a,b,c)\in R_{+}^{3}$, then $%
\phi ^{[2]}(N)(a,b,c)=N\left(
m(m(a,c),m(a,b)),m(m(a,c),m(c,b)),m(m(a,b),m(c,b))\right) $

$\leq N(m(a,c),m(a,b),m(c,b))=\phi ^{[1]}(N)(a,b,c)$. It follows by
successive iteration that $\phi ^{[n]}(N)(a,b,c)\leq \phi
^{[n-1]}(N)(a,b,c)\leq \cdots \leq N(a,b,c)$ for any positive integer $n$.
Hence $\phi ^{[n]}(N)(a,b,c)$ is a decreasing sequence of means. By Theorem 
\ref{T2}, taking the limit as $n$ approaches infinity gives $M(a,b,c)\leq
N(a,b,c)$.

Now we can prove,

\begin{theorem}
\label{cont}$M$ is continuous at each point of $R_{+}^{3}$.
\end{theorem}

\proof%
Choose continuous means $N_{1}(a,b,c)$ and $N_{2}(a,b,c)$ with $%
N_{1}(m(a,c),m(a,b),m(c,b))\leq N_{1}(a,b,c)$ and $N_{2}(a,b,c)\leq
N_{2}(m(a,c),m(a,b),m(c,b))$ for all $(a,b,c)\in R_{+}^{3}$(e.g., $%
N_{1}(a,b,c)=\min \{a,b,c\}$ and $N_{2}(a,b,c)=\max \{a,b,c\}$). Let $%
f_{n}=\phi ^{[n]}(N_{1})$ and $g_{n}=\phi ^{[n]}(N_{2})$. Then by Theorem 
\ref{T4}, $f_{n}$ is an increasing sequence of means and $g_{n}$ is a
decreasing sequence of means, each converging pointwise to $M$. Given $%
a,b,c>0$ and $\epsilon >0,$ choose $n$ so that $\left|
g_{n}(a,b,c)-f_{n}(a,b,c)\right| <\dfrac{\epsilon }{2}$. Choose $\delta >0$
so that $\left| f_{n}(a^{\prime },b^{\prime },c^{\prime
})-f_{n}(a,b,c)\right| <\dfrac{\epsilon }{2}$ and $\left| g_{n}(a^{\prime
},b^{\prime },c^{\prime })-g_{n}(a,b,c)\right| <\dfrac{\epsilon }{2}$
whenver $\left| a^{\prime }-a\right| <\delta ,\;\left| b^{\prime }-b\right|
<\delta ,\;$and $\left| c^{\prime }-c\right| <\delta $, with $a^{\prime
},b^{\prime },c^{\prime }>0$. Then $M(a^{\prime },b^{\prime },c^{\prime
})\leq g_{n}(a^{\prime },b^{\prime },c^{\prime })\leq g_{n}(a,b,c)+\dfrac{%
\epsilon }{2}<f_{n}(a,b,c)+\epsilon \leq M(a,b,c)+\epsilon $ and $%
M(a^{\prime },b^{\prime },c^{\prime })\geq f_{n}(a^{\prime },b^{\prime
},c^{\prime })\geq f_{n}(a,b,c)-\dfrac{\epsilon }{2}\geq
g_{n}(a,b,c)-\epsilon \geq M(a,b,c)-\epsilon $, which implies that $\left|
M(a^{\prime },b^{\prime },c^{\prime })-M(a,b,c)\right| <\epsilon $ whenever $%
\left| a^{\prime }-a\right| <\delta ,\;\left| b^{\prime }-b\right| <\delta ,$
and $\left| c^{\prime }-c\right| <\delta $, with $a^{\prime },b^{\prime
},c^{\prime }>0$. 
\endproof%

\begin{theorem}
\label{T3}Suppose that $m(a,b)$ is analytic in $R_{+}^{2}$. Then $M(a,b,c)$
is analytic at $Q=(s,s,s)$ for any $s>0$.
\end{theorem}

\proof%
We have that $m(z,w)$ is analytic in a region $D\subset C^{2},$ with $%
R_{+}^{2}\subset D$. Let $s>0$ be given. For $r>0,$ let $P_{r}$ equal the
polydisk $\{(z,w):\left| z-s\right| <r\}\times \{(z,w):\left| w-s\right|
<r\} $, and let Let $B_{r}$ equal the Euclidean ball $\{(z,w):$ $\left\|
(z-s,w-s)\right\| _{2}<r\}$. Choose $r>0$ suffciently small so that $m$ is
analytic at each point of $P_{r}$ and $B_{r}$. Also, for $r>0$ suffciently
small $\left| m_{z}\right| <$ $\lambda $ and $\left| m_{w}\right| <$ $%
\lambda $ at all points of $B_{r}$, with $\lambda =\dfrac{1}{\sqrt{2}}$.
This is possible since $m_{z}(s,s)=m_{w}(s,s)=\dfrac{1}{2}$. Let $a$ and $b$
be any points in $P_{r},$ and let $L$ be the line segment in $C^{2}$
connecting $(s,s)$ to $(a,b).$ $L$ is given parametrically by $z=at+s(1-t),$ 
$w=bt+s(1-t),$ $0\leq t\leq 1$. Let $g(t)=m(L(t))\Rightarrow
g(0)=m(L(0))=m(s,s)=s$. Now $g^{\prime
}(t)=m_{z}(L(t))(a-s)+m_{w}(L(t))(b-s)\left| g^{\prime }(t)\right| \leq $ $%
\lambda \sqrt{(a-s)^{2}+(b-s)^{2}}=\lambda \left\| (a-s,b-s)\right\| _{2},$ $%
0\leq t\leq 1$. Since $g(1)-s=\int_{0}^{1}g^{\prime }(t)dt,$

$\left| g(1)-s\right| \leq \lambda \left\| (a-s,b-s)\right\| _{2}$. Since $%
g(1)-s=m(a,b)-s,$ we have

\[
(a,b)\in B_{r}\Rightarrow \left| m(a,b)-s\right| \leq \lambda \left\|
(a-s,b-s)\right\| _{2}
\]

Now if $(a,b)\in P_{r}$, then $\left\| (a-s,b-s)\right\| _{2}\leq \sqrt{2}%
\left\| (a-s,b-s)\right\| _{\infty }<\sqrt{2}r\Rightarrow \left|
m(a,b)-s\right| <\lambda \sqrt{2}r=r$. Hence it follows that given any three
points $a,b,c,$ with $(a,b),$ $(a,c),$ and $(b,c)$ each in $P_{r},$ $%
(m(a,b),m(a,c)),$ $(m(a,b),m(b,c)),$ and $(m(a,c),m(b,c))$ are also each in $%
P_{r}$. Now for any positive integer $n$ and any mean $N$ analytic in $P_{r}$%
, the sequence of functions $f_{n}(a,b,c)=\phi ^{[n]}(N)(a,b,c)$ is also
analytic in $P_{r}$. We have also just shown that $\{f_{n}\}$ is uniformly
bounded on $P_{r}$. Thus $\{f_{n}\}$ has a subsequence which converges
uniformly on compact subsets of $P_{r}$ to a function $f$ analytic in $P_{r}$%
(see \cite{HL}). Since uniform convergence on compact subsets implies
pointwise convergence, Theorem \ref{T2} implies that $f(a,b,c)=M(a,b,c)$ at
all points of $P_{r}$. Therefore $M$ is analytic in $P_{r},$ and thus is
analytic at $Q$ for any $s>0$.%
\endproof%

\begin{remark}
Using the proof above, one can actually show that $M$ is analytic at any
point $(s\pm \epsilon _{1},s\pm \epsilon _{2},s\pm \epsilon _{3})$ for $%
\epsilon _{1},\epsilon _{2},$ and $\epsilon _{3}$ sufficiently small and
depending on $s$. One can then enlarge the set of points where $M$ is
analytic using the invariance. However, we have not been able to prove that $%
M$ is analytic at \textit{all} points of $R_{+}^{3}$. This is very likely
true, but a proof along the lines above has not worked so far.
\end{remark}

\section{The Invariant Logarithmic Mean}

Let $L_{3}(a,b,c)$ denote the invariant mean for the logarithmic mean $%
L(a,b)=\dfrac{b-a}{\log b-\log a}$. Note that $L$ is strict and isotone(see,
for example, \cite{CA} ), so that the results of \S \ref{s2} apply with $%
m(a,b)=L(a,b)$. The invariance property makes $L_{3}$ in some ways a natural
generalization of the logarithmic mean in two variables. By Theorem \ref{T3}
and Proposition \ref{P2}, $L_{3}$ is a strict, isotone, homogeneous mean
which is analytic at $(a,a,a)$ for any $a>0$. Using the iteration (\ref{1}),
one can fairly easily compute $L_{3}(a,b,c)$ for any $(a,b,c)\in R_{+}^{3}$.
It is well known(see, for example, \cite{CA} ) that $G_{2}(a,b)\leq
L(a,b)\leq A_{2}(a,b)$. Using Proposition \ref{P2}, it follows immediately
that $G_{3}(a,b,c)\leq L_{3}(a,b,c)\leq A_{3}(a,b,c)$. Without this
inequality, $L_{3}$ would not be a reasonable generalization of $L$.
However, we can obtain tighter upper bounds by considering the means $%
A_{p}(a,b)=\left( \dfrac{a^{p}+b^{p}}{2}\right) ^{1/p}$ and $%
A_{p}(a,b,c)=\left( \dfrac{a^{p}+b^{p}+c^{p}}{3}\right) ^{1/p}$. It has been
shown(see \cite{LI}) that $L(a,b)\leq A_{1/3}(a,b)$. Since, for any given $p$%
, $A_{p}(a,b,c)$ is invariant for $A_{p}(a,b)$, by Proposition \ref{P2}
again, it follows that $L_{3}(a,b,c)\leq A_{1/3}(a,b,c)=\left( \dfrac{%
a^{1/3}+b^{1/3}+c^{1/3}}{3}\right) ^{3}$, which is a much better bound than $%
\dfrac{a+b+c}{3}$. For example, $A_{1/3}(1,2,3)\approx \allowbreak
1.\,87934\,$, while $L_{3}(1,2,3)\approx 1.87917$

Stolarsky has defined two generalizations of $L(a,b)$ using second order
divided differences(see \cite{ST}), $U_{0}(a,b,c)=$\newline
$\left( \dfrac{1}{2}\allowbreak \dfrac{\left( a-c\right) \left( -c+b\right)
\left( -b+a\right) }{a\ln b-a\ln c+\left( \ln a\right) c-\left( \ln b\right)
c-\left( \ln a\right) b+\left( \ln c\right) b}\right) ^{1/2}$\newline
and $U_{1}(a,b,c)=\dfrac{1}{2}\dfrac{(b-c)(a-c)(a-b)}{a(b-c)\log
a-b(a-c)\log b+c(a-b)\log c}$.

$U_{1}$ is also a special case of logarithmic means in $n$ variables defined
by Pittenger(see \cite{PT}), as well as a special case of a family of means
defined by the author(see \cite{H}). It is unlikely that $U_{0}$ or $U_{1}$
are invariant means(in the next section we show that \textit{not} \textit{all%
} means in three variables are invariant). However, there is strong evidence
for

\begin{conjecture}
\label{conj1}$U_{0}(a,b,c)\leq L_{3}(a,b,c)\leq U_{1}(a,b,c)$ for all $%
(a,b,c)\in R_{+}^{3}$, with equality if and only if $a=b=c$.
\end{conjecture}

To prove Conjecture \ref{conj1}, by Theorem \ref{T4}, it suffices to prove

\begin{conjecture}
\label{conj2} $U_{0}(L(a,c),L(a,b),L(c,b))\geq U_{0}(a,b,c)$ and \newline
$U_{1}(L(a,c),L(a,b),L(c,b))\leq U_{1}(a,b,c)$ for all $(a,b,c)\in R_{+}^{3}$%
.
\end{conjecture}

There is strong numerical evidence that Conjecture \ref{conj2} is true, but
our various attempts at proving it have failed so far.

It is interesting to note that natural generalizations to three variables of
certain means $m(a,b)$ are not always comparable to the corresponding
invariant mean. For example, consider the Lehmer means $lh(a,b)=\dfrac{%
a^{2}+b^{2}}{a+b}$ and $LH(a,b,c)=\dfrac{a^{2}+b^{2}+c^{2}}{a+b+c}.$ Note
that $LH$ is not invariant for $lh$. If we let $M$ denote the invariant mean
for $lh$, then $M(1,2,3)<LH(1,2,3)$, but $M(1,2,4)>LH(1,2,4)$. Hence $M$ and 
$LH$ are not comparable.

\section{Going in Reverse\label{reverse}}

In this section we are given a mean $M(a,b,c)$, and we want to discuss the
existence and properties of a mean $m(a,b)$ such that $M$ is invariant for $%
m $. Recall that if such an $m$ exists, we say that $M$ is an invariant
mean. It is natural to ask whether there are any means in three variables
which are \textit{not} invariant. We give a simple example shortly of a mean
in three variables which is not invariant. First we discuss the analogs of
many of the results of \S \ref{s2}. We proved in \S \ref{s2} that $M$
inherits many of the properties of $m$, such as isotonicity and homogeneity.
We now prove that $m$ also inherits many of the properties of $M$, at least
when $m$ is analytic. First we need some results about means in general.

We now prove the analog of Corollary \ref{C1}.

\begin{theorem}
\label{unique}Let $M(a,b,c)$ be a differentiable, symmetric mean, and let $%
m_{1}(a,b)$ and $m_{2}(a,b)$ be analytic means, with $M$ invariant for $m_{1}
$ and for $m_{2}$. Then $m_{1}=m_{2}$.
\end{theorem}

\proof%
In general, if $M$ is invariant for $m$, then by letting $c=a$(or $%
c\rightarrow a$)$,$ we have the relation 
\begin{equation}
M(a,m(a,b),m(a,b))=M(a,a,b)  \label{5}
\end{equation}
, which holds for all $a,b>0$. For fixed $a>0$, let $f(b)=m(a,b),$ \newline
$P_{1}=(a,m(a,b),m(a,b))$, $P_{2}=(a,a,b)$, and $Q=(a,a,a)$.

By successively differentiating (\ref{5}) $k$ times, one gets an equation of
the form $(M_{x}(P_{1})+M_{y}(P_{1}))f^{(k)}(b)+$ terms which involve lower
order derivatives of $f$.

For example, differentiating (\ref{3}) twice\footnote{%
One derivative gives no information} with respect to $b$ gives

\begin{equation}
(M_{y}+M_{z})(P_{1})f^{\prime \prime
}(b)+(M_{yy}+2M_{yz}+M_{zz})(P_{1})(f^{\prime }(b))^{2}=M_{zz}(P_{2})
\label{5.5}
\end{equation}
Letting $b\rightarrow a$ and using Lemma \ref{sym}, Theorem \ref{m}(i), and
Proposition \ref{0.1}((i) and (ii)) yields $\dfrac{2}{3}f^{\prime \prime
}(a)+\frac{1}{4}M_{xx}(Q)=M_{xx}(Q)$, and hence 
\begin{equation}
f^{\prime \prime }(a)=\frac{9}{8}M_{xx}(Q)  \label{7}
\end{equation}

Differentiating (\ref{5.5}) again with respect to $b$ gives

\[
(M_{y}+M_{z})(P_{1})f^{\prime \prime \prime
}(b)+3(M_{yy}+2M_{yz}+M_{zz})(P_{1})f^{\prime }(b)f^{\prime \prime }(b)+
\]

\begin{equation}
(M_{yyy}+3M_{yyz}+3M_{yzz}+M_{zzz})(P_{1})(f^{\prime }(b))^{3}=M_{zzz}(P_{2})
\label{7.5}
\end{equation}

Again, letting $b\rightarrow a$ and using Lemma \ref{sym} and Proposition 
\ref{0.1}((i) and (ii)) yields

\begin{equation}
\dfrac{2}{3}f^{\prime \prime \prime }(a)+\dfrac{3}{2}M_{xx}(Q)f^{\prime
\prime }(a)-\dfrac{3}{4}M_{xxx}(Q)+\dfrac{3}{4}M_{xxy}(Q)=0  \label{8}
\end{equation}

It is clear that we can do this for each $k$, since we get an equation which
can be solved uniquely for $f^{(k)}(a),\;k=1,2,3,...$ If $m$ is analytic,
then $f(b)=m(a,b)$ is analytic, and thus this defines $m(a,b)$ uniquely for
each fixed $a$. 
\endproof%

\begin{theorem}
Let $M(a,b,c)$ be a symmetric mean, and let $m(a,b)$ an analytic mean, with $%
M$ invariant for $m$. If $M$ is homogeneous, then $m$ is homogeneous.
\end{theorem}

\proof%
Since $M$ is homogeneous, $M(ka,kb,kc)=kM(a,b,c)$ for any constant $k>0$. By
the invariance property (\ref{inv3}), $%
M(m(ka,kc),m(ka,kb),m(kb,kc))=kM(m(a,c),m(a,b),m(b,c))$, which implies,
again by the homogeneity of $M$, that $M(\frac{1}{k}m(ka,kc),\frac{1}{k}%
m(ka,kb),\frac{1}{k}m(kb,kc))=M(m(a,c),m(a,b),m(b,c))$. Thus $M$ is
invariant for the means $m_{1}(a,b)=m(a,b)$ and $m_{2}(a,b)=\dfrac{1}{k}%
m(ka,kb)$. By Theorem \ref{unique}, $m_{1}=m_{2}$, and thus $%
m(ka,kb)=km(a,b) $, which means that $m$ is homogeneous.

\subsection{Example}

Consider the mean $M(a,b,c)=\left( \dfrac{ab+ac+bc}{3}\right) ^{1/2}$. We
shall show that $M$ is not an invariant mean. So let $m(a,b)$ be any mean
such that $M$ is invariant for $m$. Then $M(m(a,b),m(a,b),b)=M(a,b,b)$ for
any $b>0,b\neq a$. This implies that $\left( \dfrac{m^{2}(a,b)+2bm(a,b)}{3}%
\right) ^{1/2}=\left( \dfrac{b^{2}+2ab}{3}\right) ^{1/2}\Rightarrow
m^{2}(a,b)+2bm(a,b)=b^{2}+2ab\Rightarrow $

$m(a,b)=-b\pm \sqrt{2(b^{2}+ab)}\Rightarrow M(m(1,2),m(1,3),m(2,3))=$%
\linebreak $\dfrac{1}{3}\sqrt{\left( 63-30\sqrt{3}\sqrt{2}-36\sqrt{3}+36%
\sqrt{2}-15\sqrt{3}\sqrt{2}\sqrt{5}+18\sqrt{2}\sqrt{5}+36\sqrt{5}\right) }%
\approx \allowbreak 1.\,9245$, while $M(1,2,3)=\allowbreak \dfrac{1}{3}\sqrt{%
11}\sqrt{3}\approx \allowbreak 1.\,9149$. Thus there is no mean $m(a,b)$
such that $M$ is invariant for $m$. Note that we did not need to assume that 
$m$ was analytic.

\section{Series Expansion}

Let $M(a,b,c)$ be invariant with respect to $m(a,b),$ where we assume that $M
$ is differentiable, and both symmetric and homogeneous. We would like to
expand $M$ in a Taylor Series about $Q=(1,1,1)$. Letting $a=1$ and $%
f(b)=m(1,b),$ by (\ref{7}) and Proposition \ref{0.1}(ii) 
\begin{equation}
M_{xx}(Q)=\dfrac{8}{9}f^{\prime \prime }(1),M_{xy}(Q)=-\dfrac{4}{9}f^{\prime
\prime }(1)  \label{9}
\end{equation}
The other second order partials follow from Lemma \ref{sym}. To compute the
third order partials, by Theorem \ref{m}(ii), Proposition \ref{0.1}(iii),
Proposition \ref{0.2}(i), and (\ref{9}) we have

\[
M_{xxx}(Q)=\dfrac{32}{27}\left( (f^{\prime \prime }(1))^{2}-f^{\prime \prime
}(1)\right) ,\;M_{xxy}(Q)=-\dfrac{4}{27}\left( 4(f^{\prime \prime
}(1))^{2}-f^{\prime \prime }(1)\right) 
\]

\begin{equation}
M_{xyz}(Q)=\dfrac{4}{27}\left( 8(f^{\prime \prime }(1))^{2}+f^{\prime \prime
}(1)\right)   \label{10}
\end{equation}

Finally we wish to compute the fourth order partials of $M$. Toward this
end, by Proposition \ref{0.1} and Proposition \ref{0.2}(with $a=1$), it is
easy to show that 
\begin{equation}
M_{xxxx}(Q)=-2M_{xxxy}(Q)-2M_{xxx}(Q)\text{ and }%
M_{xxyy}(Q)=-M_{xxxy}(Q)-2M_{xxy}(Q)+\frac{2}{3}M_{xyz}(Q)  \label{10.5}
\end{equation}
Differentiating (\ref{7.5}) with respect to $b$ gives

\[
(M_{y}+M_{z})(P_{1})f^{\prime \prime \prime \prime
}(b)+4(M_{yy}+2M_{yz}+M_{zz})(P_{1})f^{\prime }(b)f^{\prime \prime \prime
}(b)+
\]

\[
3(M_{yy}+2M_{yz}+M_{zz})(P_{1})(f^{\prime \prime
}(b))^{2}+6(M_{yyy}+3M_{yyz}+3M_{yzz}+M_{zzz})(P_{1})(f^{\prime
}(b))^{2}f^{\prime \prime }(b)
\]

\begin{equation}
+(M_{yyyy}+4M_{yyyz}+6M_{yyzz}+4M_{yzzz}+M_{zzzz})(P_{1})(f^{\prime
}(b))^{4}=M_{zzzz}(P_{2})  \label{10.6}
\end{equation}

By Lemma \ref{sym}, Theorem \ref{m}, Proposition \ref{0.1}, and Proposition 
\ref{0.2}, (\ref{9}), (\ref{10}), and (\ref{10.5}), (\ref{10.6}) becomes

\begin{equation}
M_{xxxy}(Q)=-\frac{16}{45}f^{(iv)}(1)-\dfrac{64}{35}\left( f^{\prime \prime
}(1)\right) ^{3}+\dfrac{448}{405}\left( f^{\prime \prime }(1)\right) ^{2}+%
\dfrac{464}{405}f^{\prime \prime }(1)  \label{10.7}
\end{equation}

By (\ref{10}), and (\ref{10.5}) we also have 
\[
M_{xxxx}(Q)=\frac{32}{45}f^{(iv)}(1)+\dfrac{128}{35}\left( f^{\prime \prime
}(1)\right) ^{3}-\dfrac{1856}{405}\left( f^{\prime \prime }(1)\right) ^{2}+%
\dfrac{32}{405}f^{\prime \prime }(1)
\]

\[
M_{xxyy}(Q)=\dfrac{16}{45}f^{(iv)}(1)+\dfrac{64}{135}\left( f^{\prime \prime
}(1)\right) ^{3}+\frac{352}{405}\left( f^{\prime \prime }(1)\right) ^{2}-%
\frac{544}{405}f^{\prime \prime }(1)
\]

\[
M_{xxyz}(Q)=-\dfrac{64}{81}\left( f^{\prime \prime }(1)\right) ^{2}-\dfrac{8%
}{81}f^{\prime \prime }(1)
\]

Again, the other third and fourth order partials follow from Lemma \ref{sym}.

Using the above, one can easily compute the Taylor polynomials of orders $%
2,3,4$, expanded about $(1,1,1)$, for the Invariant Logarithmic Mean $%
L_{3}(a,b,c)$. For example, $T_{2}(x,y,z)=\dfrac{1}{3}(x+y+z)-\dfrac{2}{27}%
\left( (x-1)^{2}+(y-1)^{2}+(z-1)^{2}\right) +$\newline
$\dfrac{2}{27}\left( (x-1)(y-1)+(x-1)(z-1)+(y-1)(z-1)\right) $.

$T_{3}(x,y,z)=T_{2}(x,y,z)+\dfrac{28}{729}\left(
(x-1)^{3}+(y-1)^{3}+(z-1)^{3}\right) $\newline
$-\dfrac{5}{243}\left(
(x-1)(y-1)(x+y-2)+(x-1)(z-1)(x+z-2)+(y-1)(z-1)(y+z-2)\right) $\newline
$+\dfrac{2}{243}(x-1)(y-1)(z-1)$

Note that $T_{3}(.9,1,1.1)=\allowbreak .\,99777\,77777\,77778\approx
\allowbreak .\,99778$, while $L_{3}(.9,1,1.1)\approx .997771$. However,
unless one is fairly close to $(1,1,1),$ $T_{3}$ and $T_{4}$ do not seem to
give that good an estimate. For example, $T_{3}(1,2,3)\approx \allowbreak
2.0 $ and $T_{4}(1,2,3)\approx \allowbreak 1.\,71251\,42073\,2902$, while $%
L_{3}(1,2,3)\approx 1.8792;$ To get a more accurate estimate for $%
L_{3}(1,2,3)$, it is better to compute $T_{4}(.5,1,1.5)$ and then use the
homogeneity property. This gives $L_{3}(1,2,3)=2L_{3}(.5,1,1.5)\approx
2T_{4}(.5,1,1.5)\approx $

$\allowbreak 1.\,88014\,30531\,0602$

\section{Another Type of Invariance}

In this section we consider the following invariance property for \textit{%
symmetric} means in two and three variables.

\begin{definition}
$M(a,b,c)$ is said to be type 2 invariant with respect to $m(a,b)$ if 
\begin{equation}
M(a,b,m(a,b))=m(a,b)  \label{inv2}
\end{equation}
for all $(a,b)\in R_{+}^{2}$.
\end{definition}

Let us call the type of invariance defined earlier as Type 1 invariance. We
shall now use the following notation: Write 
\[
(M,m)_{j}\text{ if }M(a,b,c)\text{ is type }j\text{ invariant with respect
to }m(a,b),j=1,2.
\]

It is easy to see that (\ref{inv2}) holds for the arithmetic and geometric
means in two and three variables. Indeed, let $h(u)$ be any function
monotonic on $(0,\infty )$, and define $\overline{m}(a,b)=h^{-1}m(h(a),h(b))$%
, $\overline{M}(a,b,c)=h^{-1}M(h(a),h(b),h(c))$. It follows that if $%
(M,m)_{2},$ then $(\overline{M},$ $\overline{m})_{2}$. In particular, $%
h^{-1}\left( \dfrac{h(a)+h(b)+h(c)}{2}\right) $is type2 invariant with
respect to $h^{-1}\left( \dfrac{h(a)+h(b)}{2}\right) $ for each monotone $h$.

Earlier we noted that if $lh(a,b)=\dfrac{a^{2}+b^{2}}{a+b}$ and $LH(a,b,c)=%
\dfrac{a^{2}+b^{2}+c^{2}}{a+b+c}$, then $LH$ is not type 1 invariant with
respect to $lh$. However, it is easy to show that $LH$ is type 2 invariant
with respect to $lh$. It is natural to ask whether type 1 invariance is
stronger than type 2 invariance, i.e., does type 1 invariance imply type 2
invariance ? The answer can be seen by looking at the invariant logarithmic
mean $L_{3}(a,b,c)$ defined earlier. If $L(a,b)$ denotes the logarithmic
mean, then $L_{3}(1,2,L(1,2))\approx 1.442708,$ while $L(1,2)\approx
1.442695 $. While $L_{3}(1,2,L(1,2))$ and $L(1,2)$ are close, it appears
that they are not equal. We shall prove this shortly using a series
expansion, and thus $L_{3}$ is not type 1 invariant with respect to $L$.

\begin{remark}
Unlike type 1 invariance, $m$ symmetric does not necessarily imply that $M$
is symmetric.
\end{remark}

\begin{remark}
As with type 1 invariance, one can attempt to view the invariant mean $M$ as
a special case of compounding three means in three variables by letting $%
M_{1}(a,b,c)=a,$ $M_{2}(a,b,c)=b,$and $M_{3}(a,b,c)=m(a,b)$. However, the
limit in this case does not exist.
\end{remark}

\subsection{Given $M(a,b,c)$}

We now prove an \textit{existence} result for type 2 invariance, given $%
M(a,b,c)$.

\begin{theorem}
\label{e2}Let $M(a,b,c)$ be an isotone, symmetric mean. Then there exists a
symmetric mean $m(a,b)$ such that $(M,m)_{2}$.
\end{theorem}

\proof%
Let $n(a,b)$ be any symmetric mean, and let $a$ and $b$ be given positive
numbers. Let $g(z)=M(a,b,z),$ and define the recursive sequence

$c_{k+1}=g(c_{k}),$ $c_{0}=n(a,b).$ Since $M$ is isotone, a simple inductive
proof shows that $\{c_{k}\}$ is decreasing if $c_{1}\leq c_{0},$ while $%
\{c_{k}\}$ is increasing if $c_{1}\geq c_{0}.$ Since $\left| c_{k}\right|
\leq \max \{a,b\},$ in either case $\{c_{k}\}$ is bounded and monotonic, and
hence converges to some real number $L,$ $\min \{a,b\}\leq L\leq b\leq \max
\{a,b\}.$ Note that if $c_{1}=c_{0},$ then $c_{k}=c_{0}$ for all $k,$ and
thus $\{c_{k}\}$ converges to $L=c_{0}$. Define the mean $m(a,b)=L$. Since $%
M $ and $n$ are symmetric, $m(b,a)=m(a,b)\Rightarrow m$ is symmetric. Of
course the iteration $c_{k+1}=g(c_{k})$ converges to a fixed point of $g$,
and thus $g(L)=L.$ This implies that $M(a,b,m(a,b))=m(a,b),$ which proves
that $(M,m)_{2}$. The only thing left to prove is that $m$ is continuous.
Toward this end, let $\{n_{k}\}$ be the sequence of means defined by the
recursion $n_{k+1}(a,b)=M(a,b,n_{k}(a,b)),$ $n_{0}(a,b)=n(a,b).$ Note that
since $M$ and $n$ are continuous, each $n_{k}$ is also continuous. It is
easy to show that, for each fixed $a$ and $b,$ $c_{k+1}=n_{k}(a,b)$. Now let 
$g_{k}(a,b)$ be the sequence of means corresponding to $n(a,b)=\min \{a,b\},$
and let $h_{k}(a,b)$ be the sequence of means corresponding to $n(a,b)=\max
\{a,b\}.$ If $c_{0}=\min \{a,b\},$ then $c_{1}\leq c_{0}$ and thus $%
\{c_{k}\} $ is decreasing. Hence $g_{k}$ is a decreasing sequence of means
converging to $m(a,b).$ Similarly, $h_{k}$ is an increasing sequence of
means converging to $m(a,b).$ The rest of the proof now follows in a similar
fashion to the proof of Theorem \ref{cont}. Given $a,b>0$ and $\epsilon >0,$
choose $k$ so that $\left| g_{k}(a,b)-f_{k}(a,b)\right| <\dfrac{\epsilon }{2}
$. Choose $\delta >0$ so that $\left| g_{k}(a^{\prime },b^{\prime
})-g_{k}(a,b)\right| <\dfrac{\epsilon }{2}$ and $\left| h_{k}(a^{\prime
},b^{\prime })-h_{k}(a,b)\right| <\dfrac{\epsilon }{2}$ whenever $\left|
a^{\prime }-a\right| <\delta ,\;$and $\left| b^{\prime }-b\right| <\delta $,
with $a^{\prime },b^{\prime }>0$. Then $m(a^{\prime },b^{\prime })\leq
g_{k}(a^{\prime },b^{\prime })\leq g_{k}(a,b)+\dfrac{\epsilon }{2}%
<h_{k}(a,b)+\epsilon \leq m(a,b)+\epsilon $ and

$m(a^{\prime },b^{\prime })\geq h_{k}(a^{\prime },b^{\prime })\geq
h_{k}(a,b)-\dfrac{\epsilon }{2}\geq g_{k}(a,b)-\epsilon \geq m(a,b)-\epsilon 
$ , which implies that $\left| m(a^{\prime },b^{\prime })-m(a,b)\right|
<\epsilon $ whenver $\left| a^{\prime }-a\right| <\delta $ and $\;\left|
b^{\prime }-b\right| <\delta $, with $a^{\prime },b^{\prime }>0$. 
\endproof%

One may of course discuss type 2 invariance for functions $M(a,b,c)$ and $%
m(a,b)$ which are not necessarily means. Our next result shows, however,
that if $M$ is a mean, then $m$ must also be a mean.

\begin{lemma}
\label{mean}Suppose that $M(a,b,c)$ is a strict mean, and that $g(a,b)$ is
any symmetric continuous function satisfying $M(a,b,g(a,b))=g(a,b)$ for all $%
(a,b)\in R_{+}^{2}$. Then $g$ is a strict symmetric mean.
\end{lemma}

\proof%
Suppose that $a<b$ and that $g(a,b)\leq a.$ Since $M$ is a strict mean, $%
M(a,b,g(a,b))>g(a,b)=M(a,b,g(a,b)),$ a contradiction. Hence $g(a,b)>a$
whenever $a<b.$ Similarly, $g(a,b)<b$ whenever $a<b,$ and thus $g$ is a
strict mean.%
\endproof%

\begin{lemma}
\label{iso}Suppose that $M(a,b,c)$ is a strictly isotone mean, and that $%
m(a,b)$ is a mean with $(M,m)_{2}$. Then $m$ is a strictly isotone mean.
\end{lemma}

\proof%
Let $a>0$ and suppose that $b_{1}<b_{2}$ with $m(a,b_{1})=m(a,b_{2}).$ Since 
$(M,m)_{2}$, this implies that $M(a,b_{1},m(a,b_{1}))=M(a,b_{2},m(a,b_{1})),$
which contradicts the fact that $M$ is strictly isotone. Hence $m(a,b)$ is
either increasing or decreasing in $b.$ Now for $a\leq b,a=m(a,a)$ and $%
a\leq m(a,b)$. Thus $m(a,b)$ must be increasing in $b$ for any fixed $a>0.$
Similarly, $m(a,b)$ must be increasing in $a$ for any fixed $b>0.$

We have not been able to prove a \textit{uniqueness} result for type 2
invariance. That is, does $(M,m_{1})_{2}=(M,m_{2})_{2}$ imply that $%
m_{1}=m_{2}$ ? We can prove uniqueness with the additional assumption that $%
M_{zz}$ does not change sign.

\begin{remark}
\label{local}The proofs given above show that Lemmas \ref{iso}, \ref{mean},
and Theorem \ref{uniq} are all \textit{local} results. That is, one need
only assume that the hypotheses hold for all $(a,b)\in I_{1}\times I_{2}$,
where $I_{1}$ and $I_{2}$ are open subintervals of the positive reals. The
conclusions then also hold for all $(a,b)\in I_{1}\times I_{2}$.
\end{remark}

\begin{proposition}
\label{diff}Let $M(a,b,c)$ be a symmetric, strictly isotone mean which is
twice differentiable in $R_{+}^{3}$. Assume also that $M_{zz}(x,y,z)$ is
never $0$ on $R_{+}^{3}$, and let $m(a,b)$ be a mean with $(M,m)_{2}.$ Then $%
f(b)=m(a,b)$ is differentiable for each fixed $a>0,$ and $M_{z}(a,b,m(a,b))<1
$ for each $a,b>0$.
\end{proposition}

\proof%
$M_{y}(a,b,m(a,b))\geq 0$ and $M_{z}(a,b,m(a,b))\geq 0$ since $M$ is
strictly isotone. Note that $M_{yy}(a,b,c)=M_{zz}(a,c,b)$ since $M$ is
symmetric. Hence $M_{yy}(a,b,c)$ is never $0$ on $R_{+}^{3}.$ This easily
implies the strict positivity of $M_{y}$--i.e., $M_{y}(a,b,c)>0$ for all $%
(a,b,c)\in R_{+}^{3}$. By Theorem \ref{e2}, there exists a symmetric mean $%
m(a,b)$ such that $(M,m)_{2}.$ Let $f(b)=m(a,b)$ for each fixed $a>0$. Since 
$f$ is increasing by Lemma \ref{iso}, $f^{\prime }(b)$ exists on a set $%
S\subset (0,\infty ),$ with $m(S^{c})=0,$ where $S^{c}$ denotes the
complement of $S$. Differentiating both sides of (\ref{inv2}) with respect
to $b$ gives

\begin{equation}
M_{y}(a,b,m(a,b))+M_{z}(a,b,m(a,b))f^{\prime }(b)=f^{\prime }(b)  \label{db}
\end{equation}

Note that (\ref{db}) implies that $f^{\prime }(b)>0$ on $S$. (\ref{db}) also
implies that, for $b\in S,$ $M_{y}(a,b,f(b))+M_{z}(a,b,f(b))f^{\prime
}(b)>M_{z}(a,b,f(b))f^{\prime }(b)\Rightarrow $

$f^{\prime }(b)>M_{z}(a,b,f(b))f^{\prime }(b)$. Hence we have $%
M_{z}(a,b,f(b))<1$ whenever $b\in S$. Solving (\ref{db}) for $f^{\prime }(b)$
yields $f^{\prime }(b)=\dfrac{M_{y}(a,b,f(b))}{1-M_{z}(a,b,f(b))}.$ This
shows that $f^{\prime \prime }(b)$ exists(indeed $f^{(k)}(b)$ exists for any 
$k)$. Hence, if $b\in S,$ we can differentiate both sides of (\ref{db}) with
respect to $b$ to obtain

\begin{equation}
M_{yy}(P)+2M_{yz}(P)f^{\prime }(b)+M_{z}(P)f^{\prime \prime
}(b)+M_{zz}(P)(f^{\prime }(b))^{2}=f^{\prime \prime }(b)  \label{12}
\end{equation}

where $P=(a,b,f(b)).$ Thus $M_{yy}(P)+2M_{yz}(P)f^{\prime
}(b)+M_{zz}(P)(f^{\prime }(b))^{2}=(1-M_{zz}(P))f^{\prime \prime }(b).$ Upon
dividing thru by $f^{\prime }(b)$ we have 
\begin{equation}
\frac{M_{yy}(P)}{f^{\prime }(b)}+2M_{yz}(P)+M_{zz}(P)f^{\prime
}(b)=(1-M_{z}(P))\frac{f^{\prime \prime }(b)}{f^{\prime }(b)}  \label{12.5}
\end{equation}
Now suppose that $\{b_{n}\}$ is a sequence in $S$, with $b_{n}\rightarrow
b^{+}$, $b\notin S$. If $f^{\prime }(b_{n})\rightarrow r^{+},$ then by (\ref
{db}), $r=M_{y}(a,b,f(b))+M_{z}(a,b,f(b))r>M_{z}(a,b,f(b))r\Rightarrow
M_{z}(a,b,f(b))<1.$The same conclusion follows if $b_{n}\rightarrow b^{-}$.
By taking convergent subsequences, we can now conclude that

\[
\text{If }b_{n}\in S,b\notin S,b_{n}\rightarrow b^{+}\text{ or }b^{-},\text{
and }f^{\prime }(b_{n})\nrightarrow \infty \text{, then }M_{z}(a,b,f(b))<1
\]
We now prove that $f^{\prime }(b_{n})$ \textit{cannot} approach $\infty .$
So suppose that $b_{n}\in S,b_{n}\rightarrow b_{0}$ and $f^{\prime
}(b_{n})\rightarrow \infty .$ Assume first that $M_{zz}(a,b,c)>0$ on $%
R_{+}^{3}$. Then, upon replacing $b$ by $b_{n},$ the LHS of (\ref{12.5})
approaches $\infty $ as $n\rightarrow \infty .$ Thus the LHS of (\ref{12.5})
is positive for $n$ sufficiently large. Since $M_{z}(a,b_{n},f(b_{n}))<1$,
the RHS of (\ref{12.5}) implies that $f^{\prime \prime }(b_{n})>0$ for $n$
sufficiently large. If $M_{zz}(a,b,c)<0$ on $R_{+}^{3},$ then a similar
argument shows that $f^{\prime \prime }(b_{n})<0$ for $n$ sufficiently
large. Now if $f^{\prime }(b_{n})\rightarrow \infty $ as $b_{n}\rightarrow
b_{0}$ from \textit{both} sides, then $f$ must be concave to one side of $%
b_{0}$ and convex on the other side of $b_{0}$. We have just shown that that
is impossible, and thus $f^{\prime }(b_{n})\nrightarrow \infty $ as $%
b_{n}\rightarrow b_{0}.$ We can now conclude that

\begin{equation}
M_{z}(a,b,m(a,b))<1\text{ for all }a,b>0  \label{lt1}
\end{equation}
Finally, let $c=m(a,b),\vartriangle c=m(a,b+\vartriangle b)-m(a,b).$ By the
Mean Value Theorem, $M(a,b+\vartriangle b,c+\vartriangle c)-M(a,b,c)=$

$M_{y}(a,b+t\vartriangle b,c+t\vartriangle c)\vartriangle
b+M_{z}(a,b+t\vartriangle b,c+t\vartriangle c)\vartriangle c,$ $0<t<1.$
Also, since $(M,m)_{2},$ $M(a,b+\vartriangle b,c+\vartriangle
c)-M(a,b,c)=M(a,b+\vartriangle b,m(a,b+\vartriangle b)-M(a,b,m(a,b))=$

$m(a,b+\vartriangle b)-m(a,b)=\vartriangle c$. Hence $\dfrac{\vartriangle c}{%
\vartriangle b}=M_{y}(a,b+t\vartriangle b,c+t\vartriangle
c)+M_{z}(a,b+t\vartriangle b,c+t\vartriangle c)\dfrac{\vartriangle c}{%
\vartriangle b}\Rightarrow $

$\dfrac{\vartriangle c}{\vartriangle b}=\dfrac{M_{y}(a,b+t\vartriangle
b,c+t\vartriangle c)}{1-M_{z}(a,b+t\vartriangle b,c+t\vartriangle c)}.$
Letting $\vartriangle b\rightarrow 0,$ (\ref{lt1}) implies that $f^{\prime
}(b)$ exists and equals $\dfrac{M_{y}(a,b,c)}{1-M_{z}(a,b,c)}$ whenever $%
M_{z}(a,b,m(a,b))\neq 1.$%
\endproof%

\begin{theorem}
\label{uniq}Let $M(a,b,c)$ be a symmetric, strictly isotone mean which is
twice differentiable in $R_{+}^{3}$. Assume that $M_{zz}(x,y,z)$ is never $0$
on $R_{+}^{3}.$ If $m_{1}$ and $m_{2}$ are means with $%
(M,m_{1})_{2}=(M,m_{2})_{2}$ , then $m_{1}=m_{2}.$
\end{theorem}

\proof%
Let $g(z)=M(a,b,z)$ for fixed $a>0,$ $b>0$. If $M(a,b,m_{1}(a,b))=m_{1}(a,b)$
and $M(a,b,m_{2}(a,b))=m_{2}(a,b),$ then $g(L_{1})=L_{1}$ and $%
g(L_{2})=L_{2},$ where $L_{j}=m_{j}(a,b).$ If $m_{1}(a,b)\neq m_{2}(a,b),$
then $L_{1}$ and $L_{2}$ are \textit{distinct} fixed points of $g$. Now $%
\dfrac{g(L_{2})-g(L_{1})}{L_{2}-L_{1}}=1\Rightarrow g^{\prime }(t)=1$ for
some $t,$ $L_{1}<t<L_{2}$. By Proposition \ref{diff}, $g^{\prime
}(L_{1})\leq 1$ and $g^{\prime }(L_{2})\leq 1.$ Since $g$ cannot be
constant, this implies that $g^{\prime }$ must be increasing somewhere on $%
(L_{1},L_{2})$ and decreasing somewhere on $(L_{1},L_{2}).$ That contradicts
the fact that $g^{\prime \prime }$ is never $0$ on $(0,\infty ).$ Hence $%
m_{1}(a,b)=m_{2}(a,b)$ for all $a>0,$ $b>0$.%
\endproof%

We now prove the existence of an analytic $m$ such that $(M,m)_{2}.$

\begin{theorem}
\label{anal2}Let $M(a,b,c)$ be a symmetric, strictly isotone mean which is
analytic in $R_{+}^{3}$. Assume also that $M_{zz}(x,y,z)$ is never $0$ on $%
R_{+}^{3}.$ Then there exists a unique symmetric mean $m(a,b)$ which is
analytic in $R_{+}^{2},$ and such that $(M,m)_{2}.$
\end{theorem}

\proof%
Let $T(z_{1},z_{2},z_{3})=M(z_{1},z_{2},z_{3})-z_{3}$, which is analytic in
some open set in $C^{3}$ containing $R_{+}^{3}$. Then, for all $x>0,$ $y>0,$ 
$T(x,y,m(x,y))=0$ and $T_{z_{3}}(x,y,m(x,y))\neq 0$ by Proposition \ref{diff}%
. By the Implicit Function Theorem(see \cite{HL}), for any given $a>0,$ $b>0$%
, the equation $T(z_{1},z_{2},z_{3})=0$ along with $T(a,b,m(a,b))=0$ has a
unique solution $z_{3}=g(z_{1},z_{2})$ analytic in some open neighborhood $O$
of $(a,b)$ in $C^{2}$, with $g(a,b)=m(a,b)$. Restricting $z_{1}$ and $z_{2}$
to be real, we have $M(x,y,g(x,y))=g(x,y)$ for all $(x,y)\in I=O\cap
R_{+}^{2}$. Note that since $%
M(z_{2},z_{1},z_{3})-z_{3}=M(z_{1},z_{2},z_{3})-z_{3}$ and $m(b,a)=m(a,b)$,
by uniqueness $g(z_{2},z_{1})=g(z_{1},z_{2})$. By Lemma \ref{mean}(see also
Remark \ref{local}), $g(x,y)$ must be a symmetric mean, at least for $%
(x,y)\in I$. Since $M(x,y,m(x,y))=m(x,y)$ for all $(x,y)\in R_{+}^{2}$,
Theorem \ref{uniq}(see also Remark \ref{local}) then implies that $%
m(x,y)=g(x,y)$ for all $(x,y)\in I$. Hence $m$ extends to be analytic in an
open neighborhood of any $(a,b)\in R_{+}^{2}.$ The uniqueness of $m$ follows
from Theorem \ref{uniq}. 
\endproof%

\QTP{Body Math}
\begin{example}
Let $M(a,b,c)=\left( \dfrac{ab+ac+bc}{3}\right) ^{1/2}$, which interpolates
the arithmetic and geometric means in three variables. Then

\QTP{Body Math}
$M_{zz}=\allowbreak -\dfrac{1}{4}\dfrac{\left( x+y\right) ^{2}}{\left(
xy+xz+yz\right) \sqrt{\left( 3xy+3xz+3yz\right) }}$, which is never $0$ on $%
R_{+}^{3}.$ Hence, by Theorem \ref{anal2}, there exists a unique symmetric
mean $m(a,b)$ which is analytic in $R_{+}^{2},$ and such that $(M,m)_{2}.$
Let $g(a,b)=\left( \dfrac{a^{1/2}+b^{1/2}}{2}\right) ^{2}$, which equals the
Holder mean $A_{1/2}(a,b)$. Let $h(a,b)=M(a,b,g(a,b))=\allowbreak \dfrac{1}{6%
}\sqrt{\left( 18ab+3a^{2}+6\left( \sqrt{a}\right) ^{3}\sqrt{b}+6\left( \sqrt{%
b}\right) ^{3}\sqrt{a}+3b^{2}\right) }$.

\QTP{Body Math}
Now $h(1,b)=\allowbreak \dfrac{1}{6}\sqrt{\left( 18b+3+6\sqrt{b}+6\left( 
\sqrt{b}\right) ^{3}+3b^{2}\right) }$ and $g(1,b)=\dfrac{1}{4}\left( 1+\sqrt{%
b}\right) ^{2}$. Then $h(1,b)>g(1,b)\Leftrightarrow $

\QTP{Body Math}
$\dfrac{1}{36}\left( 3u^{4}+6u^{3}+18u^{2}+6u+3\right) >\dfrac{1}{16}%
(1+u)^{4}$, where $u=\sqrt{b}$. The latter inequality holds if and only if $%
(u-1)^{4}>0$, which holds for any $u\neq 1$. Hence $h(1,b)>g(1,b)$ for any $%
b\neq 1$. Since $g$ and $h$ are each homogeneous of degree $1$, $%
h(a,b)>g(a,b)$ for all $a,b>0$ with $a\neq b$. We have shown that

\QTP{Body Math}
$M(a,b,g(a,b))\geq g(a,b),$ with equality if and only if $a=b$. It follows
that if $\{n_{k}\}$ is the sequence of means defined by the recursion $%
n_{k+1}(a,b)=M(a,b,n_{k}(a,b)),$ $n_{0}(a,b)=g(a,b)$, then $\{n_{k}(a,b)\}$
is increasing and converges to $m(a,b)$ for each $a,b>0$(see the proof of
Theorem \ref{e2}). Since $n_{1}(a,b)=h(a,b)>g(a,b)=$ $n_{0}(a,b)$ for all $%
a,b>0$ with $a\neq b,$ we have proven that $m(a,b)>\left( \dfrac{%
a^{1/2}+b^{1/2}}{2}\right) ^{2},$ with equality if and only if $a=b$.
\end{example}

We now prove a uniqueness result without most of the assumptions on $M$ in
Theorem \ref{uniq}. However, we then must assume that $m_{1}$ and $m_{2}$
are analytic.

\begin{theorem}
Let $M(a,b,c)$ be a differentiable, symmetric mean, and let $m_{1}(a,b)$ and 
$m_{2}(a,b)$ be symmetric means, each analytic in $R_{+}^{2}$. Then $%
(M,m_{1})_{2}=(M,m_{2})_{2}$ implies that $m_{1}=m_{2}.$
\end{theorem}

\proof%
If $M(a,b,c)$ is type 2 invariant with respect to $m(a,b),$ then 
\begin{equation}
M(a,b,m(a,b))=m(a,b)  \label{11}
\end{equation}

Letting $f(b)=m(a,b)$ and differentiating both sides of (\ref{12}) with
respect to $b$ gives

\[
M_{z}(P)f^{\prime \prime \prime }(b)+3M_{yz}(P)f^{\prime \prime
}(b)+3M_{zz}(P)f^{\prime }(b)f^{\prime \prime }(b)+3M_{yyz}(P)f^{\prime }(b)
\]

\begin{equation}
+3M_{yzz}(P)(f^{\prime }(b))^{2}+M_{yyy}(P)+M_{zzz}(P)(f^{\prime
}(b))^{3}=f^{\prime \prime \prime }(b)  \label{14}
\end{equation}

In general, taking $k$ derivatives of (\ref{11}) with respect to $b$ yields
an equation of the form $M_{z}(P)f^{(k)}(b)+$ $A=f^{(k)}(b),$ where $A$ is a
polynomial in the partial derivatives of $M$(evaluated at $P$) and
derivatives of $f$(evaluated at $b$) of order $<k$. Since $M_{z}(a,a,a)=%
\dfrac{1}{3}$, letting $b\rightarrow a,$we get an equation which can be
solved uniquely for $f^{(k)}(a),\;k=1,2,3,...$ If $m$ is analytic, then $%
f(b)=m(a,b)$ is analytic, and thus this defines $m(a,b)$ uniquely for each
fixed $a$. 
\endproof%

\subsection{Series Expansion}

Letting $b\rightarrow a$ in (\ref{12}), and using $f^{\prime }(a)=\dfrac{1}{2%
}$ and $M_{z}(Q)=\dfrac{1}{3}$ yields $M_{yy}(Q)+M_{yz}(Q)+\dfrac{1}{4}%
M_{zz}=\dfrac{2}{3}f^{\prime \prime }(a)$, where $Q=(a,a,a)$. Since $M$ is
symmetric, Lemma \ref{sym} implies $\dfrac{5}{4}M_{xx}(Q)+M_{xy}(Q)=\dfrac{2%
}{3}f^{\prime \prime }(a)$. Finally, applying Proposition \ref{0.1} yields 
\begin{equation}
M_{xx}(Q)=\dfrac{8}{9}f^{\prime \prime }(a)  \label{13}
\end{equation}
$,$ which is just (\ref{8}).

Letting $b\rightarrow a$ in (\ref{14}) and using Lemma \ref{sym},
Proposition \ref{0.1}(ii), and Proposition \ref{0.2}(i) yields $-\dfrac{9}{8a%
}M_{xx}(Q)=\dfrac{2}{3}f^{\prime \prime \prime }(a)$. By (\ref{13}) we have $%
f^{\prime \prime \prime }(a)=-\dfrac{3}{2a}f^{\prime \prime }(a),$ which
already holds by Theorem \ref{m}. So to get new information we must
differentiate both sides of (\ref{14}) again with respect to $b.$ Then
letting $b\rightarrow a$ and using Lemma \ref{sym} gives

\[
4M_{yz}(Q)f^{\prime \prime \prime }(a)+2M_{zz}(Q)f^{\prime \prime \prime
}(a)+3M_{zz}(Q)(f^{\prime \prime }(a))^{2}+12M_{yyz}(Q)f^{\prime \prime }(a)
\]

\begin{equation}
+\dfrac{3}{2}M_{zzz}(P)f^{\prime \prime }(a)+\dfrac{3}{2}M_{yyzz}(Q)+\dfrac{5%
}{2}M_{yyyz}(Q)+\dfrac{17}{16}M_{xxxx}(Q)=\dfrac{2}{3}f^{(iv)}(a)  \label{15}
\end{equation}

Since we are interested in the series expanded about $(1,1,1)$, we now let $%
a=1$. By (\ref{10.5}), Theorem \ref{m}, and (\ref{13}), (\ref{15}) becomes 
\begin{equation}
\dfrac{8}{3}\left( f^{\prime \prime }(1)\right) ^{3}-\dfrac{16}{3}\left(
f^{\prime \prime }(1)\right) ^{2}+\dfrac{8}{3}f^{\prime \prime }(1)+\frac{3}{%
8}M_{xxx}(Q)-\frac{9}{2}f^{\prime \prime }(1)M_{xxx}(Q)-\frac{9}{8}%
M_{xxxy}(Q)=\frac{2}{3}f^{(iv)}(1)  \label{16}
\end{equation}

\subsection{Given $m(a,b)$}

We now prove a result similar to Theorem \ref{e2}, except here we are given
the mean $m(a,b).$

\begin{theorem}
Let $m(a,b)$ be a symmetric mean. Then there exists a symmetric mean $%
M(a,b,c)$ such that $(M,m)_{2}.$
\end{theorem}

\proof%
Let $n(a,b)$ be any mean in two variables(not neceassarily symmetric), and
let $M(a,b,c)=n(m(a,b),c),$ for $c$ between $a$ and $b$. Then extend $%
M(a,b,c)$ to the rest or $R_{+}^{3}$ so that $M$ is symmetric. It follows
that

\begin{equation}
M(a,b,c)=\left\{ 
\begin{array}{lll}
n(m(a,b),c) & \text{if} & a\leq c\leq b\text{ or }b\leq c\leq a \\ 
n(m(a,c),b) & \text{if} & a\leq b\leq c\text{ or }c\leq b\leq a \\ 
n(m(c,b),a) & \text{if} & b\leq a\leq c\text{ or }c\leq a\leq b
\end{array}
\right.  \label{symm}
\end{equation}

Since $m$ is symmetric, it is not hard to see that $M$ must also be
symmetric. It also follows easily that $M$ is continuous. For example, fix $%
a\neq c$, and let $b\rightarrow c$. Using the definition of $M$ above(first
two rows), $n(m(a,b),c)$ and $n(m(a,c),b)$ each approach $n(m(a,c),c)$ by
the continuity of $m$ and of $n$, respectively. Finally, since $m$ is a
mean, if $c=m(a,b),$ then $a\leq c\leq b$ or $b\leq c\leq a$. Hence $%
M(a,b,m(a,b))=n(m(a,b),m(a,b))=m(a,b).$ Thus $(M,m)_{2}$.

\begin{remark}
If $m(a,b)$ is \textit{analytic}, we have not been able to prove the
existence of a \textit{symmetric} \textit{analytic }mean $M(a,b,c)$ such
that $(M,m)_{2}.$ The Theorem above takes care of the symmetric part, and
letting $M(a,b,c)=n(m(a,b),c)$ for \textit{all }$(a,b,c)\in R_{+}^{3},$ with 
$n$ analytic, forces $M$ to be analytic. However, then $M$ is \textit{not}
symmetric in general. It is satisfying both conditions that does not seem
easy to do. Now in certain special cases it is clear how to choose $n$ so
that $m$ is both symmetric and analytic\textit{. }However, these choices do
not seem to generalize. For example, if $m(a,b)=\dfrac{a+b}{2}$, then one
can choose $n(a,b)=\dfrac{2}{3}a+\dfrac{1}{3}b$, which implies that $%
n(m(a,b),c)=\allowbreak \dfrac{1}{3}a+\dfrac{1}{3}b+\dfrac{1}{3}c,$ which of
course is the arithmetic mean in three variables. Similarly, if $m(a,b)=%
\sqrt{ab}$, then one can choose $n(a,b)=a^{2/3}b^{1/3}$, which implies that $%
n(m(a,b),c)=\allowbreak \sqrt[3]{abc}.$ The appearance of the $\dfrac{2}{3}$
and $\dfrac{1}{3}$ in each case is not an accident. It is easy to show that
a necessary condition for the $M$ from (\ref{symm}) above to be
differentiable is $n_{y}(a,a)=\dfrac{1}{2}n_{x}(a,a)$. Also, note that in
each of these cases, $n(a,b)=(m\otimes _{a}m)(b,a)$, where $m\otimes _{a}n$
denotes the Archimedean compound of $m$ with $n$(see \cite{B}). One might be
tempted to try this approach in general, but it does not work. For example,
if $m(a,b)=\dfrac{a^{2}+b^{2}}{a+b}$ and $n(a,b)=(m\otimes _{a}m)(b,a)$,
then $n(m(1,2),3)\approx 2.307961$

and $n(m(1,3),2)\approx 2.356566$. Hence if $M(a,b,c)=n(m(a,b),c)$ for all%
\textit{\ }$(a,b,c)\in R_{+}^{3}$, then $M$ is not symmetric. Of course
there \textit{is} a mean $M(a,b,c)$ such that $(M,m)_{2}$--namely $M(a,b,c)=%
\dfrac{a^{2}+b^{2}+c^{2}}{a+b+c}$.
\end{remark}

\subsection{Invariance of Both Types}

It is also interesting to consider means $M(a,b,c)$ and $m(a,b)$ with $M$ 
\textbf{both} type 1 and type 2 invariant with respect to $m$. By (\ref{10}%
), $M_{xxx}(Q)=\dfrac{1}{27}\left( 32(f^{\prime \prime }(1))^{2}-32f^{\prime
\prime }(1)\right) ,$ $Q=(1,1,1).$ Substituting into (\ref{16}) yields 
\begin{equation}
M_{xxxy}(Q)=-\frac{8}{9}\left( \frac{2}{3}f^{(iv)}(1)+\dfrac{8}{3}\left(
f^{\prime \prime }(1)\right) ^{3}-\dfrac{4}{9}\left( f^{\prime \prime
}(1)\right) ^{2}-\dfrac{20}{9}f^{\prime \prime }(1)\right)  \label{17}
\end{equation}

Now let $m(a,b)=L(a,b)=\dfrac{b-a}{\ln b-\ln a}$, $f(x)=\dfrac{x-1}{\ln x}$

By (\ref{17}), if $M$ is \textbf{both} type 1 and type 2 invariant with
respect to $m(a,b),$ then

$M_{xxxy}(Q)=-\dfrac{8}{9}\lim\limits_{x\rightarrow 1}\left( \dfrac{2}{3}%
f^{(iv)}(x)+\dfrac{8}{3}\left( f^{\prime \prime }(x)\right) ^{3}-\dfrac{4}{9}%
\left( f^{\prime \prime }(x)\right) ^{2}-\dfrac{20}{9}f^{\prime \prime
}(x)\right) =\allowbreak \dfrac{248}{3645}$

$\approx \allowbreak 6.\,80381495\times 10^{-2}$

By (\ref{10.7}), if $M(a,b,c)=L_{3}(a,b,c),$ the invariant logarithmic mean,
which is type 1 invariant with respect to $m(x,y)$, then

$M_{xxxy}(Q)=\dfrac{8}{15}\lim\limits_{x\rightarrow 1}\left( -\dfrac{2}{3}%
f^{(iv)}(x)-\dfrac{8}{9}\left( f^{\prime \prime }(x)\right) ^{3}+\dfrac{56}{%
27}\left( f^{\prime \prime }(x)\right) ^{2}+\dfrac{58}{27}f^{\prime \prime
}(x)\right) =\allowbreak \dfrac{136}{2025}$

$\approx \allowbreak 6.\,71604\,9\times 10^{-2}$

Since the two values are not equal, $L_{3}$ is \textbf{not} type 2 invariant
with respect to $L(x,y)=\dfrac{x-y}{\ln x-\ln y}$. In fact, we believe the
following is true.

\begin{conjecture}
Let $m(a,b)=\dfrac{a+b}{2}$ and $M(a,b,c)=\dfrac{a+b+c}{3}.$ Then the 
\textit{only} pairs of means $\{\overline{m}(a,b),\overline{M}(a,b,c)\}$
which satisfy \textbf{both} $(\overline{M},\overline{m})_{1}$ and $(%
\overline{M},\overline{m})_{2}$ are the means $\overline{m}%
(a,b)=h^{-1}m(h(a),h(b))$, $\overline{M}(a,b,c)=h^{-1}M(h(a),h(b),h(c)),$
where $h(u)$ is a function monotonic on $(0,\infty ).$
\end{conjecture}

\section{Open Questions and Future Reserach}

(1) One can, of course, define invariance for \textbf{nonsymmetric} means.
For example, if $m(a,b)=\dfrac{2}{3}a+\dfrac{1}{3}b$ and $M(a,b,c)=\dfrac{4}{%
7}a+\dfrac{2}{7}b+\dfrac{1}{7}c$, then

$M(a,b,m(a,b))=\allowbreak \dfrac{2}{3}a+\dfrac{1}{3}b=m(a,b)$ and

\QTP{Body Math}
$M(m(a,b),m(a,c),m(b,c))=\allowbreak \dfrac{4}{7}a+\dfrac{2}{7}b+\dfrac{1}{7}%
c=M(a,b,c)$. Hence $(M,m)_{1}$ and $(M,m)_{2}$. Note, however, that

\QTP{Body Math}
$M(m(a,c),m(a,b),m(b,c))=\allowbreak \dfrac{4}{7}a+\dfrac{5}{21}c+\dfrac{4}{%
21}b\neq m(a,b).$

\QTP{Body Math}
(2) One might discuss invariance of types 1 and 2 for classes of functions
other than means. For example, polynomials or rational functions. One can
explore questions such as: What are the invariant polynomials or rational
functions in three variables ?

\QTP{Body Math}
For Type 1 Invariance:

\QTP{Body Math}
(3) Given any symmetric mean $m(a,b),$ is there always a mean $M(a,b,c)$
such that $(M,m)_{1}$ ? We proved this with the additional assumptions that $%
m$ is strict and isotone.

\QTP{Body Math}
(4) As discussed earlier, show that $U_{0}(a,b,c)\leq L_{3}(a,b,c)\leq
U_{1}(a,b,c)$ for all $(a,b,c)\in R_{+}^{3},$ where $L_{3}$ is the invariant
logarithmic mean and $U_{0}$ and $U_{1}$ are Stolarsky's generalizations of
the logarithmic mean $L(a,b).$

\QTP{Body Math}
For Type 2 Invariance:

\QTP{Body Math}
(5) As discussed earlier, given a symmetric analytic mean $m(a,b),$ is there
always a symmetric analytic mean $M(a,b,c)$ such that $(M,m)_{2}$ ? In
particular, is there an analytic mean $M(a,b,c)$ which is Type 2 invariant
with respect to the logarithmic mean $L(a,b)$.

\QTP{Body Math}
(6) Given a symmetric isotone mean $M(a,b,c),$ and means $m_{1}(a,b)$ and $%
m_{2}(a,b)$ with $(M,m_{1})_{2}$ and $(M,m_{2})_{2},$ must $m_{1}=m_{2}$ ?
We proved this with additional assumptions on $M$ and/or $m$.

\end{document}